\font\math=msbm10 scaled 1200

\newcommand{\cpd}{{\hbox{\math o}}}

\documentstyle[12pt]{article}
\font\germ=eufm10 at12pt
\def\goth#1{\hbox{\germ#1}}
\newcommand{\sdag}{\scriptsize \dag}

\newcommand{\qed}{\hbox{\rule[-2pt]{3pt}{6pt}}}

\newtheorem{Thm}{Theorem}[section]

\newtheorem{rem}{Remark}[section]

\newtheorem{ex}{Example}[section]

\newtheorem{defi}{Definition}[section]
\newtheorem{lem}{Lemma}[section]

\newtheorem{cor}{Corollary}[section]

\newcommand{\qedh}{\hfill\qed}

\newcommand{\nears}{\!\!\!\!\!\!\!\!\!\!\!}
\newcommand{\vv}{ \vspace{.3in}  }

\def\labelenumi{\theenumi.}
\def\theenumi{\arabic{enumi}}

\def\labelenumi{\theenumi.}
\def\theenumi{\Alph{enumi}}
\renewcommand{\theenumi}{\Alph{enumi}}

\def\labelenumi{\theenumi.}
\def\theenumi{\arabic{enumi}}

\def\labelenumi{\theenumi.}
\def\theenumi{{\rm (\roman{enumi})}}

\begin{document}








\newpage

\begin{center}

{\Large {\bf  Infinitesimal Takesaki duality \\of
 Hamiltonian
 vector fields \\ on a symplectic manifold}}

\end{center}
\vspace{.5in}
\begin{center}
{\large Katsunori Kawamura
\footnote{ e-mail : kawamura@kurims.kyoto-u.ac.jp.}}

\end{center}
\begin{center}
{\normalsize  Research Institute for Mathematical Sciences }

{\normalsize  Kyoto University, Kyoto 606,Japan}
\end{center}

\begin{center}
{\large \today .}
\end{center}
\vspace{.4in}
\begin{abstract}

For an infinitesimal 
symplectic action of a Lie algebra ${\goth g}$ 
on a symplectic manifold,
 we construct an infinitesimal crossed product of
 Hamiltonian vector fields and Lie algebra ${\goth g}$.
 We obtain its second crossed product in case
 ${\goth g}={\bf R}$ and show an 
 infinitesimal version for a 
theorem type of Takesaki duality.

\end{abstract}

\section{Introduction}

By an application of a functional representation for 
 ${\cal L}({\cal H})$ 
 \cite{CMP94},  
 any von Neumann algebra is faithfully represented as
 functions on a projective Hilbert space 
${\cal P}({\cal H})$ by $*$-product
 as an unital W$^{*}$-algebra.
 Its non commutativity comes from the symplectic structure
 of ${\cal P}({\cal H})$ induced by the fundamental form of
 a  K\"{a}hler
 metric of ${\cal P}({\cal H})$.
 These results are showed in section 2.
  In this view, we expect that
 several structure of von Neumann algebra
 come from symplectic geometry.
 
 Our aim is a  study how the structures of
 von Neumann algebra is recognized as the symplectic geometry
 on ${\cal P}({\cal H})$.
 
For example, by \cite{CMP90} the $*$-product of functions
 of some K\"{a}hler manifold $M$ is depend on its 
 holomorphic sectional curvature $c$ of $M$
 by the deformation parameter $1/c$. 
 And the functional representation for  von Neumann algebra
 is the case of $c=1$. And 
 the represented commutator of von Neumann algebra
 becomes just Poisson bracket of functions multiplying $\sqrt{-1}$.
 In other words, von Neumann algebra
 is a special kind of " symplectic algebra ".
 ( See corollary \ref{cor:21}, too. )

 In this paper, 
 we define an  infinitesimal crossed product
 by the semi direct product of Lie algebras and 
 show Takesaki duality type theorem for
 some dynamical system on a symplectic manifold.

\begin{Thm}(Infinitesimal Takesaki duality)

Let $M$ be a symplectic manifold and $(M,{\bf R},\beta)$
 an infinitesimal symplectic dynamical system 
 induced by some symplectic dynamical system,  
${\goth X}_{H}(M)$ the Lie algebra of 
all Hamiltonian vector fields on $M$.
 Then an  isomorphism of Lie algebras  follows
\[ ({\goth X}_{H}(M)\cpd_{\beta}{\bf R})\cpd_{\hat{\beta}}{\bf R}
\cong {\goth X}_{H}(M)\cpd_{\delta}H_{1}\]
 where $\hat{\beta}$ is the dual action of $\beta$ and
 $H_{1}$ is the Heisenberg Lie algebra with two generators
$a, a^{\sdag}$ over ${\bf C}$  which satisfying
$[a,a^{\sdag}]=I$.  
The left hand side is semi direct product
 of Lie algebras induced by a derivative   action 
$\delta : H_{1}\to 
{\rm End}({\goth X}_{H}(M))$ such that 
$\delta_{a}X_{f}=\frac{1}{\sqrt{2}}X_{\beta_{1}f}, 
\delta_{a^{\sdag}}X_{f}=-\frac{1}{\sqrt{2}}X_{\beta_{1}f}$
 for $X_{f}\in {\goth X}_{H}(M), f\in C^{\infty}(M)$.

\end{Thm}

By the functional representation for  von Neumann algebra,
 this theorem is recognized as a generalization of 
 Takesaki duality for non associative algebra
 related to symplectic structure of symplectic manifold.
 About Takesaki duality, see \cite{Takesaki}.

\section{ A functional representation for  
non commutative von Neumann algebras }

We construct a functional representation of 
von Neumann algebras as follows.

Let ${\cal H}$ be a Hilbert space
 and ${\cal L}({\cal H})$ the set of 
all bounded linear operators on ${\cal H}$.
  The following result are obtained by 
Cirelli, Lanzavecchia, Mani\`{a}-Pizzocchero et al.

\begin{Thm} (A functional representation for  
${\cal L}({\cal H})$ \cite{CMP94})\label{Thm:21}
\begin{enumerate}
\item ${\cal P}({\cal H})\equiv({\cal H}-\{0\})/{\bf C}^{\times}$
 becomes a K\"{a}hler manifold $({\cal P}({\cal H}),g)$ by the 
inhomogeneous coordinate and the 
Fubini-study type metric $g$.
\item A set of functions 
 ${\cal K}({\cal P}({\cal H}))\equiv 
\{ f\in C^{\infty}({\cal P}({\cal H})) : D^{2}f=0, \bar{D}^{2}f=0\}$
 becomes a C$^{*}$-algebra by a $*$-product
\[ (f*l)({\bf x})\equiv f({\bf x}) l({\bf x})
+\partial_{{\bf x}} f({\rm grad}_{{\bf x}}l) \hspace{.2in}
 f,l\in {\cal K}({\cal P}({\cal H})), {\bf x}\in{\cal P}({\cal H})\]
, a $*$-involution 
  \[f^{*}\equiv \mbox{ complex conjugate } \bar{f}\mbox{ of } f\]
and a C$^{*}$-norm
\[ \|f\|\equiv\sup_{{\bf x}\in {\cal P}({\cal H})}
|(\bar{f}*f)({\bf x})|^{1/2}.\]
Here $D,\bar{D}$ are the holomorphic and anti-holomorphic part of covariant
 derivative, and ${\rm grad}f$ is the holomorphic  part of
  complexified gradient of  $f$ by the K\"{a}hler metric $g$.
\item ${\cal L}({\cal H})\cong {\cal K}({\cal P}({\cal H}))$ as
 a  C$^{*}$-algebra.
\end{enumerate}
\end{Thm}
{\it Proof.}  (i) See \cite{CLM}. 
 (  Topology is induced by inhomogeneous coordinate.
 Differential structure is defined by Frechet differential.)
 (ii) and  (iii) is 
in the section 2 of \cite{CMP94}.
 For $A\in {\cal L}({\cal H})$, let $f_{A}([x])\equiv
 <x|Ax>$ for $[x]\in {\cal P}({\cal H})$, $x\in {\cal H}$, $\|x\|=1$.
 \qedh

\begin{rem}\label{rem:21}{\rm 
 In this theorem, K\"{a}hler manifold means a Hilbert
 manifold with a nondegenerate hermitian
 metric  which fundamental form is closed.
 And if ${\rm dim}{\cal H}=n\in {\bf N}$,
 then ${\cal P}({\cal H})$
 becomes ordinary projective space ${\bf C}P^{n-1}$
 with the standard K\"{a}hler structure.
}
\end{rem}

By the definition   of Poison bracket $\{\cdot,\cdot\}$
 of the fundamental form of ${\cal P}({\cal H})$,
 the following relation holds

\begin{equation}\label
{eqn:21}
f*l-l*f=\sqrt{-1}\{f,l\}
\end{equation}
for $f,l\in {\cal K}({\cal P}({\cal H}))$. 
\begin{defi}
\begin{enumerate}
\item For a subset  ${\cal S}\subset {\cal K}({\cal P}({\cal H}))$,
  the {\it  Poisson commutant }  ${\cal S}^{c}
\equiv\{ f\in {\cal K}({\cal P}({\cal H}))
 : \{f,l\}=0 \mbox{ for } l\in {\cal S}\}$. 
\item The strong K\"{a}hler topology of ${\cal K}({\cal P}({\cal H}))$
 is defined by pointwise convergence 
 $f_{\lambda}\to f$ of net $\{f_{\lambda}\}$ of functions
 in ${\cal K}({\cal P}({\cal H}))$
 by $(\bar{f_{\lambda}}*f_{\lambda})(p)\to (\bar{f}*f)(p)$ for each $p\in 
 {\cal P}({\cal H})$.
\item The weak  K\"{a}hler topology of ${\cal K}({\cal P}({\cal H}))$
 is defined by pointwise convergence 
$f_{\lambda}\to f$ of net $\{f_{\lambda}\}$ 
 of functions by $f_{\lambda}(p)\to f(p)$ for each $p\in {\cal P}({\cal H})$.
\end{enumerate}
\end{defi}
Remark if we say " subalgebra " of ${\cal K}({\cal P}({\cal H}))$,
 it means a subalgebra by $*$-product in theorem \ref{Thm:21} (ii). 

 By using  theorem \ref{Thm:21}, we construct a 
functional faithful 
representation for 
a  non commutative von Neumann algebras.
\begin{cor}\label{cor:21}( A functional representation for 
 a  von Neumann algebra)
For  each von Neumann algebra acting ${\cal R}$ on 
 a Hilbert space ${\cal H}$, there is a subalgebra ${\cal S}$
 of ${\cal K}({\cal P}({\cal H}))$ such that ${\cal S}^{cc}={\cal S}$
 and closed under complex conjugate and it is isomorphic to ${\cal R}$
 as an unital W$^{*}$-algebra by the strong K\"{a}hler  topology in
  ${\cal K}({\cal P}({\cal H}))$.
\end{cor}
{\it Proof.}  
By theorem \ref{Thm:21}  (iii), 
for a von Neumann algebra ${\cal R}$,
 let  ${\cal S}\subset 
{\cal K}({\cal P}({\cal H}))$ be
 a range of the functional representation of 
${\cal R}\subset{\cal L}({\cal H})$. Then 
 ${\cal S}$ becomes $*$-subalgebra by theorem  \ref{Thm:21} (iii)
 and by the relation of commutator and Poisson bracket
 in the last argument, 
it satisfies  ${\cal S}^{cc}={\cal S}$.
 
Because of definition of norm
 in $ {\cal K}({\cal P}({\cal H}))$,
 $\|Ax\|=|(\bar{f_{A}}*f_{A})({\bf x})|$ for $A\in {\cal R}$
, $x\in {\cal H}, {\bf x}=[x]\in {\cal P}({\cal H})$ and
 $f_{A}$ is a function in ${\cal S}$ 
 corresponding to $A$.
The topology of ${\cal S}$ related to
 strong operator topology of ${\cal R}$
 is just equal to 
the strong K\"{a}hler topology.  \qedh

\begin{cor}\label{cor:22}
For a  $*$-subalgebra ${\cal S}$ 
with $1_{M}$ of ${\cal K}({\cal P}({\cal H}))$,
 ${\cal S}$ is  weak K\"{a}hler  closed if and only if
 ${\cal S}={\cal S}^{cc}$.
\end{cor}
 {\it Proof.}
 By corollary \ref{cor:21}, ${\cal S}={\cal S}^{cc}$ if and only if
 there is a related von Neumann algebra  ${\cal R}$ and it is 
 weak operator closed.  By definition of  
weak  K\"{a}hler  topology, it is equivalent to weak operator topology 
 ${\cal R}$. Since a unital $*$-algebra ${\cal R}$ of ${\cal L}({\cal H})$
 is weak operator closed if and only if it is strong operator closed, 
 ${\cal S}$ is weak K\"{a}hler  closed if and only if ${\cal S}$
 is  strong K\"{a}hler  closed.
  \qedh
\vv

\paragraph{ Comparison with non-commutative geometry}

In the theory of non-commutative geometry by A.Connes
 \cite{Connes1}, operator algebras are treated as like some function
 spaces on '' virtual ''  manifolds.
 there are many applications and similar the theory of results for
 the theory of original geometry \cite{Connesbook}.

But we cannot view non-commutative manifold
 which related to non-commutative operator algebras.

 Corollary \ref{cor:21} gives some answers for this problems as follows:
 By the functional representation, any von Neumann algebra
 acting on a Hilbert space ${\cal H}$ 
 is identified with function $*$-sub algebra
 ${\cal S}$ of ${\cal K}({\cal P}({\cal H}))$
 such that it is closed under bi-commutator
 in  ${\cal K}({\cal P}({\cal H}))$.
 ( About  the case of C$^{*}$-algebras, see \cite{CMP94}. ) 
 By an equation \ref{eqn:21}, the commutator of functions is
 given by Poisson bracket on ${\cal P}({\cal H})$.
 Since Poisson bracket is defined by symplectic form 
 ( fundamental form ) of ${\cal P}({\cal H})$.

 We can conclude that {\it
 '' non commutativity of von Neumann algebras
 comes from skew symmetry of symplectic form 
on ${\cal P}({\cal H})$ ''}.

symplectic structure is well known and fundamental non commutative structure of 
 non ''  virtual '' geometry. This is a solution
 for how non-commutative structure of operator algebra 
 comes from real geometry.

\begin{rem}\label{rem:22}
{\rm
\vspace{-.1in}
\begin{enumerate}
\item  Usually, von Neumann algebra is considered
 as " non commutative $L_{p}$-space " or some kind of dual of measure
 space. So someone may feel that
 measurable functions  on measure spaces are more suitable for
 a functional representation of von Neumann algebra rather
 than smooth functions on a  smooth manifold. 
But we don't touch such a way of view  in this paper. 
\vspace{-.1in}
\item For a unital C$^{*}$-algebra ${\goth A}$ 
 acting on a Hilbert space
 ${\cal H}$, its enveloping von Neumann algebra
  $\overline{{\goth A}}$ is the  strong operator closure of it.
 So someone may feel that a functional representation
 is  given by extension of Cirelli-Mani\`{a}-Pizzocchero
 functional representation
 by
 taking closure of it \cite{CMP94}. But it does not stands
 because generally 
Cirelli-Mani\`{a}-Pizzocchero functional representation
 is not strong operator continuous.
\vspace{-.1in}
\item 
The relations of so  called $*$-deformation of  functions 
 over a symplectic manifold 
and C$^{*}$-algebra are studied by Riefell \cite{Rieffel}.
 In our paper, we use a   $*$-product as similar as $*$-deformation, but 
 there are several differences :
(a) There is no deformation parameter 
 in appearance in our theory.
 In \cite{CMP90}, deformation parameter is 
equal to the  inverse of holomorphic 
sectional curvature of K\"{a}hler manifold.
 In this case deformation parameter is fixed to 1 because
 projective space has holomorphic sectional curvature 1.
 (b) We treat only first-order deformation.
Usually, the deformation algebra is the algebra
 of formal power series with coefficient of functions on 
 s symplectic manifold. But our case, 
there is no term except 0-th  and 1-st orders.
 (c) Our  deformed functional space is not
 closed usual point wise product of functions.
 (d) We don't deform operator algebras but represent  as function
 space.
  \end{enumerate}}
\end{rem}

\section{Infinitesimal symplectic dynamical system}

We define the infinitesimal symplectic dynamical system
 and introduce some example of it.
--text follows this line--
 Let $M$ be a symplectic manifold i.e. a smooth manifold with
 nondegenerate closed 2-form $\omega$ and
 the set ${\goth X}(M)$ of all
 smooth vector fields on $M$. ${\rm Symp}(M)$ is
 the group of all symplectic diffeomorphisms on $M$.
For 
$f\in C^{\infty}(M)$, a vector field $X_{f}$ is called
 Hamiltonian if it satisfies 
$\omega_{p}((X_{f})_{p},u)=d_{p}f(u)$ for any $p\in M$
 and any tangent vector $u$ at $p$. It is well defined
 by non degeneracy of $\omega$. We denote the  set of all
 Hamiltonian vector fields on $M$ by ${\goth X}_{H}(M)$.
 Then ${\goth X}_{H}(M)$ is a Lie subalgebra of ${\goth X}(M)$
 by the commutator $[\cdot,\cdot]$ of vector fields.
 and $(C^{\infty}(M),\{\cdot,\cdot\})\to 
({\goth X}_{H}(M) ,[\cdot,\cdot])$, $f\mapsto -X_{f}$ 
 is  surjective
 Lie homomorphism with kernel ${\bf C}1_{M}$
 where $\{\cdot,\cdot\}:
C^{\infty}(M)\times C^{\infty}(M)\to C^{\infty}(M)$ 
is the Poisson bracket defined by
 $\{f,l\}_{p}\equiv\omega_{p}((X_{f})_{p}, (X_{l})_{p})$ 
for $f,l\in C^{\infty}(M)$ at $p\in M$ \cite{Fom}.

 Let $G$ be a Lie  group with its Lie algebra ${\goth g}$.
\begin{defi}
\begin{enumerate}
\item $(M,G,\alpha)$ is a symplectic dynamical system
 if $\alpha:G\to {\rm Symp}(M)$ is a symplectic action of $G$ on $M$ 
 which is smooth  $G\ni g\mapsto \alpha_{g}(p)\in M$ at  each
 point $p\in M$.
\item $(M,{\goth g},\beta)$ is an infinitesimal 
symplectic dynamical system
 if $\beta : {\goth g}\to {\goth X}(M)$ is Lie homomorphism
 such that 
$\beta_{X}\{f,l\}=\{\beta_{X}f,l\}+\{f,\beta_{X}l\}$ for
 any $f,l\in C^{\infty}(M)$ and $X\in {\goth g}$.
\end{enumerate}
\end{defi}

\begin{rem} \label{rem:31}
{\rm

For a given symplectic dynamical system $(M,G,\alpha)$,
 $(\beta_{X}f)(p)\equiv \\
\frac{d}{dt}(\alpha_{{\rm exp}(tX)}^{*}f)(p)|_{t=0}$
 for $p\in M$ 
 becomes an infinitesimal 
symplectic dynamical system. The reason why we treat infinitesimal 
symplectic dynamical systems is follows :
 Ordinary crossed product
 in the theory of operator algebras  is an associative algebra generated
 by an associative algebra and a group algebra. But we treat
 $C^{\infty}(M)$ as a Lie algebra by its Poisson bracket. So we consider
 to construct a crossed product of Lie algebras as a Lie algebra
 by using semi direct product of Lie algebras.
 Relation between von Neumann algebra and infinitesimal 
symplectic dynamical system is explained  later.}
\end{rem}

\begin{ex}{\rm

(i) $G={\bf R}_{+}\equiv\{e^{t} : t\in {\bf R}\}$, ${\goth g}={\bf R}$.
 Then $\alpha$ is symplectic flow on $M$.

(ii) If $M$ is a K\"{a}hler manifold and $\alpha:G\to{\rm  Iso}_{h}(M)$
 is holomorphic K\"{a}hler isometric action, then automatically
 $\alpha$ becomes symplectic action
 for the fundamental form of
 the K\"{a}hler metric. So $(M,G,\alpha)$ becomes
 symplectic dynamical system. So, if $G$ is a Lie subgroup of
 K\"{a}hler isometries which are holomorphic gives 
 a symplectic dynamical system.

(iii) Specially simply connected, connected complete K\"{a}hler manifold
 $M$ 
 with  non zero constant holomorphic sectional curvature $c$ is considered 
  in some geometric generalization of 
  quantum mechanics 
(\cite{ACLM}, \cite{CL}, \cite{CMP90}). Its symplectic dynamical
 system $(M,{\bf R}_{+},\alpha)$ induces the 
generalized Schr\"{o}dinger equation
 on $M$. 
In this case, we need to  remove a condition of smoothness of action
 ${\bf R}_{+}\ni x\mapsto \alpha_{x}(p)$ for $p\in M$
 without some region 
 when it is necessary to treat  unboundedness of Hamiltonian.

(iv) If $M={\cal P}({\cal H})$
over ( possibly  
 infinite dimensional )  Hilbert space ${\cal H}$,
 its symplectic dynamical system is
 related to W$^{*}$-dynamical system
 $({\cal R}, G, \sigma)$
 ( appendix \ref{section:A} )
 by a functional representation of 
 a von Neumann algebra ${\cal R}$ on ${\cal P}({\cal H})$
 ( appendix \ref{section:B} )
 where $\sigma$ is an action of locally compact group
 $G$ on ${\cal R}$ as $*$-automorphisms implemented
 by unitary representation of $G$ on ${\cal H}$.

 In the case $G$ is a Lie group with the Lie algebra ${\goth g}$
 and  there is a dense linear subspace  ${\cal D}\subset {\cal H}$
 such that $dU:{\goth g}\to {\rm End}({\cal D})$ 
 is the infinitesimal representation of ${\goth g}$
 induced by the smooth representation  $U$ of $G$ on ${\cal D}$.

Let $d\sigma_{X}\equiv {\rm ad}_{dU_{X}}$ for $x\in{\goth g}$.
 Then $d\sigma:{\goth g}\to {\rm End}({\cal R})$
 becomes representation of ${\goth g}$
 where ${\cal R}$ is restricted on ${\cal D}$.
 By the functional representation ${\cal S}=\{ f_{A} : A\in {\cal R} \}$ 
of ${\cal R}$,  define
 $\beta_{X}(f_{A})\equiv f_{d\sigma_{X}(A)}$
 for $X\in {\goth g}$.
 Then $({\cal P}({\cal H}), {\goth g}, \beta)$ becomes
 an infinitesimal dynamical system with function space ${\cal S}\subset
 {\cal K}({\cal P}({\cal H}))$.
 
}
\end{ex}

\section{Infinitesimal symplectic covariant system}

We define infinitesimal symplectic covariant system and
 introduce some example, and
 define crossed product of 
 Hamiltonian vector fields.

Let $(M,{\goth g},\beta)$ be an infinitesimal symplectic dynamical system.

\begin{defi}
$(\pi, \Lambda)$ is an infinitesimal covariant system
 of $(M,{\goth g},\beta)$ if there is a vector space $V$,
 $\pi : (C^{\infty}(M),\{\cdot,\cdot\})\to {\rm End}(V), \Lambda :
 {\goth g}\to {\rm End}(V)$ are Lie homomorphisms and they satisfy
 $[\Lambda_{X}, \pi(f)]=\pi(\beta_{X}f)$ for each $X\in {\goth g}, 
 f\in C^{\infty}(M)$.
\end{defi}

Clearly, infinitesimal symplectic covariant systems
 come from symplectic covariant systems
 by differential ${\bf R}\ni t\mapsto {\rm exp}(tX)\in G$.
  
Once more we remark that the reason
 why we use infinitesimal form is in order to construct crossed
 product  of Lie algebras of $(C^{\infty}(M),\{\cdot,\cdot\})$ 
 and some  Lie algebra because von Neumann algebra
 is a sub Lie algebra of $C^{\infty}({\cal P}({\cal H}))$
 by Poisson bracket.

\begin{ex}\label{ex:41}
{\rm
\begin{enumerate}
 \item  Let \[V\equiv C^{\infty}(M)\hspace{.3in}  
\pi(f)l\equiv ({\rm ad}f)(l)=-X_{f}l=\{f,l\}\hspace{.3in}  
     \Lambda_{X}l\equiv \beta_{X}l\]
 for $X\in {\goth g}, f,l\in C^{\infty}(M)$. Then $(\pi,\Lambda)$
 becomes an infinitesimal symplectic covariant system
 by Leibniz rule of $\beta_{X}$.

\item Let $(M,G,\alpha)$ be a symplectic dynamical system
 which induces  a given  infinitesimal symplectic dynamical system
$(M,{\goth g},\beta)$. 
Let 
\[V\equiv C^{\infty}(M)\otimes C^{\infty}(G)\]
 \[(\pi(f)(l\otimes \chi))(p,g)\equiv 
\{\alpha_{g^{-1}}^{*}f,l\}_{p}\chi(g)\hspace{.3in} 
 \Lambda_{X}\equiv I\otimes l_{X}\]
 for $X\in {\goth g}, f,l\in C^{\infty}(M)$
, $\chi\in C^{\infty}(G)$
 where $\otimes$ is algebraic tensor over ${\bf C}$
 and  $(l_{X}\chi)(g)\equiv \frac{d}{dt}
\chi({\rm exp}(-tX)\cdot g)|_{t=0}$.
 Then $(\pi,\Lambda)$
 becomes an infinitesimal symplectic covariant system.
 As the  same way, $\tilde{V}\equiv C^{\infty}(M\times G)$ gives
  an infinitesimal symplectic covariant system. If
 we take the right regular representation,
 \[(\pi^{'}(f)(l\otimes \chi))(p,g)\equiv 
\{\alpha_{g}^{*}f,l\}_{p}\chi(g)\hspace{.3in} 
 \Lambda^{'}_{X}\equiv I\otimes r_{X}\]
 where  $X\in {\goth g}, f,l\in C^{\infty}(M)
, \chi\in C^{\infty}(G)$,
 $(r_{X}\chi)(g)\equiv \frac{d}{dt}
\chi( g\cdot {\rm exp}(tX))|_{t=0}$.
 Then $(\pi^{'},\Lambda^{'})$
  it becomes
  an infinitesimal symplectic covariant system, too.
\end{enumerate}}

\end{ex}

\begin{rem}\label{rem:42}{\rm
In the theory of operator algebras, a covariant system is 
 a couple of representations of algebra and group  defined
 on a same Hilbert space. But in the infinitesimal case,
 we can not use Hilbert space generally. For example,
 consider time development of unbounded Hamiltonian $H$, 
 then its differential is not defined on whole
 Hilbert space. As same as in Tomita-Takesaki theory,
 $H={\rm log}\Delta$ is unbounded operator generally
 where $\Delta$ is a modular operator.
 So generator of action of ${\bf R}$
 is not able to defined on whole of $L_{2}({\bf R})$.
}
\end{rem}

In order to  define a crossed product of 
an infinitesimal symplectic dynamical system $(M,{\goth g},\beta)$ 
by an infinitesimal  covariant system $(\pi,\Lambda)$,
 we review the semi direct product of Lie algebras.

For two Lie algebras ${\goth g}_{1},{\goth g}_{2}$,
 consider a representation   
$\delta :{\goth g}_{2}\to {\rm End}({\goth g}_{1})$
  such that 
$\delta_{b}$  is a derivation that is 
$\delta_{b}([v,w])=[\delta_{b}v,w]+[v, \delta_{b}w]$
 for $b\in{\goth g}_{2}, v,w\in {\goth g}_{1}$.
 We call such $\delta$ by derivative action of ${\goth g}_{2}$
 on ${\goth g}_{1}$.
 
 The semi direct product 
 ${\goth g}_{1}\cpd_{\delta}{\goth g}_{2}$
 is a Lie algebra which is equal to 
direct sum ${\goth g}_{1}\oplus {\goth g}_{2}$
 as a linear space and has a  new Lie bracket
 $[b, v]\equiv \delta_{b}v$ 
for $b\in{\goth g}_{2}, v\in {\goth g}_{1}$ and other Lie brackets
 are  same as that of ${\goth g}_{1},{\goth g}_{2}$ respectively.

For 
 an infinitesimal symplectic dynamical system 
$(M,{\goth g},\beta)$, $\beta :{\goth g}\to {\goth X}(M)\subset
{\rm End}(C^{\infty}(M))$ is a derivative action.
For an infinitesimal  covariant system $(\pi,\Lambda)$
 of $(M,{\goth g},\beta)$, 
if let $\gamma_{\xi}\equiv{\rm ad}\xi$
 for $\xi \in\Lambda({\goth g})$, 
 then
 $\gamma:\Lambda({\goth g})\to {\rm End}(\pi(C^{\infty}(M))$
 is a derivative action, too. So, we can construct
 two semi direct products
 $C^{\infty}(M)  \cpd_{\beta}{\goth g}$
 and $\pi(C^{\infty}(M))
\cpd_{\gamma} \Lambda({\goth g})$.

If $(\pi,\Lambda)$ is the infinitesimal  covariant system
 in previous example (ii), then
$\pi(C^{\infty}(M))\cong  {\goth X}_{H}(M)$ because
 kernel of $\pi={\bf C}1_{M}$ and 
 $\Lambda({\goth g})\cong {\goth g}$. Therefore
\begin{defi}(Crossed product of Hamiltonian vector field)
${\goth X}_{H}(M)  \cpd_{\beta}{\goth g}
 \equiv  
 \pi(C^{\infty}(M))
\cpd_{\gamma} \Lambda({\goth g})$
 for the infinitesimal  covariant system
 in example \ref{ex:41} (ii).
\end{defi}
%



\section{Second crossed product}

To state the duality theorem,
 we need
the definition of second crossed product and 
  restrict $G$  to  an abelian Lie group 
 which dual becomes Lie algebras,too because
 the duality theorem comes from
 comes from duality  of abelian group and we treat   a smooth class of
 group. Then we treat only $G={\bf R}_{+}$.
Let $(M,{\bf R},\beta)$ be an infinitesimal  dynamical system.
 For ${\cal L}_{1}\equiv {\goth X}_{H}(M)\cpd_{\beta}{\bf R}$, 
we define a
 dual infinitesimal covariant  system. Remark 
 now $G={\bf R}_{+}, {\goth g}={\bf R}$. But
 we identify $G$ and ${\goth g}$ in this case from now.

 Let $V_{1}\equiv C^{\infty}(M)\otimes C^{\infty}({\bf R})$
 and define an  operator $W_{s} : V_{1}\to V_{1}$ by 
\[(W_{s} \Psi)(p,t)\equiv {\rm exp}(\sqrt{-1}ts)
\Psi (p,t)\]  
for $\Psi\in V_{1}, s,t\in {\bf R}$
 and 
\[\hat{\alpha}_{s}\equiv {\rm Ad}W_{s}\hspace{.4in}
  \hat{\beta}_{s}\equiv \frac{d}{dx}\hat{\alpha}_{sx}|_{x=0}\]
 on ${\cal L}_{1}$.

Let $(w_{s}\chi)(t)\equiv 
\sqrt{-1}st\chi(t)$ for 
 $\chi \in C^{\infty}({\bf R})$ and $t,s\in{\bf R}$. 
Then  $((I\otimes w_{s})\Psi)(p,t)= 
  \frac{d}{dx}(W_{s} \Psi)(p,t)|_{x=0}$,
  and  $\hat{\beta}_{s}(R)=[I\otimes w_{s},R]$ for $s\in{\bf R}, 
R\in {\cal L}_{1}$. Furthermore 
 by Jacobi identity of bracket, 
$\hat{\beta} :{\bf R}\to {\rm End}({\cal L}_{1})$ 
 becomes  a derivative action.

\begin{lem}
\begin{enumerate}
\item $\hat{\beta}_{s}(\pi(f))=0$ for $f\in C^{\infty}(M),  s\in {\bf R}$.
\item $\hat{\beta}_{s}(\Lambda_{t})=\sqrt{-1}stI$
\end{enumerate}
\end{lem}
{\it Proof}
(i) By definition, $\hat{\alpha}_{s}(\pi(f))=\pi(f)$ 
for any $f\in C^{\infty}(M), s\in{\bf R}$. So, the statement follows.
 
(ii) $(\hat{\alpha}_{s}(\Lambda_{t})(\Psi))(p,r)$

\[\begin{array}{ll}
=&
(W_{s}\Lambda_{t}W_{-s}\Psi)(p,r) \\
=&e^{\sqrt{-1}rs}(\Lambda_{t}W_{-s}\Psi)(p,r)\\
=&e^{\sqrt{-1}rs}\frac{d}{dx}(W_{-s}\Psi)(p,r-tx)|_{x=0}\\
=&e^{\sqrt{-1}rs}\frac{d}{dx}e^{-\sqrt{-1}(r-tx)s}\Psi(p,r-tx)|_{x=0}\\
=&\sqrt{-1}ts\Psi(p,r)+\frac{d}{dx}\Psi(p,r-tx)|_{x=0}.\\
\end{array}\]
So $\hat{\alpha}_{s}(\Lambda_{t})=\Lambda_{t}+\sqrt{-1}stI$.
Therefore , the statement follows.\qedh

By this lemma, $\hat{\alpha}_{t}([R,R^{'}])=[R,R^{'}]$
 for each $R,R^{'}\in {\cal L}_{1}, t\in {\bf R}$.

Let $V_{2}\equiv V_{1}\otimes C^{\infty}({\bf R})$.
Define $\hat{\pi} : {\cal L}_{1}
\to {\rm End}(V)$
 and $\hat{\Lambda}:{\bf R}\to {\rm End}(V)$
 by $(\hat{\pi} (R)\Phi)(p,s,r)
\equiv(\hat{\alpha}_{-r}(R)\otimes I)\Phi)(p,s,r)$
 and $\hat{\Lambda}_{t}=I\otimes I\otimes l_{t}$
 for $r,s,t\in {\bf R},R\in {\cal L}_{1}, \Phi\in V_{2}$.

\begin{lem}
\begin{enumerate}
\item $[\hat{\Lambda}_{t},\hat{\pi}(R)]
=\hat{\pi}(\hat{\beta}_{t}(R))$
 on $V_{2}$ for $t\in {\bf R}, R\in {\cal L}_{1}$
\item $\hat{\pi}(\pi(f))=\pi(f)\otimes I$ for $f\in C^{\infty}(M)$
\item $\hat{\pi}(\Lambda_{t})=\Lambda_{t}\otimes I-I\otimes I \otimes w_{t}$
 for $t\in {\bf R}$ for $t\in {\bf R}$.
 
\end{enumerate}
\end{lem}
{\it Proof}
(i) 
For  $\Phi\in V_{2}$,

$([\hat{\Lambda}_{t},\hat{\pi}(R)]\Phi)(p,r,s)$
\[
\begin{array}{lll}
&= & \frac{d}{dx}(\hat{\pi}(R)\Phi)(p,r,s-xt)|_{x=0}\\

&&- \{(\hat{\alpha}_{-s}(R)\otimes I)\hat{\Lambda}_{t}\Phi\}(p,r,s) \\

&=& \frac{d}{dx}\{(\hat{\alpha}_{-(s-xt)}(R)\otimes I)\Phi\}(p,r,s-xt)|_{x=0}\\
&&- \{(\hat{\alpha}_{-s}(R)\otimes I)\frac{d}{dx}\Phi\}(p,r,s-xt)|_{x=0}\\

&= &\{(\hat{\alpha}_{-s}(\frac{d}{dx}
\hat{\alpha}_{xt}(R))|_{x=0}\otimes I)\Phi\}(p,r,s)\\

&=&\{\hat{\alpha}_{-s} (\hat{\beta}_{t}(R))\otimes I)\Phi\}(p,r,s)\\
&=& \{\hat{\pi}(\hat{\beta}_{t}(R))\Phi\}(p,r,s)
 
\end{array}.\]

(ii) By proof of lemma 5.1 (i), it follows.

(iii) $(\hat{\pi}(\Lambda_{t})\Phi)(p,r,s)
=((\hat{\alpha}_{-s}(\Lambda_{t})\otimes I)\Phi)(p,r,s)$

By proof of lemma 5.1 (ii),

$=\{((\Lambda_{t}-\sqrt{-1}st)\otimes I)\Phi\}(p,r,s)$

$=((\Lambda_{t}\otimes I-I\otimes I\otimes w_{t})\Phi)(p,r,s)$
\qedh

We define second crossed product 
$({\goth X}_{H}(M)\cpd_{\beta}{\bf R})\cpd_{\hat{\beta}}{\bf R}$
 by semi direct product 
$\hat{\pi}({\goth X}_{H}(M)\cpd_{\beta}{\bf R}  )\cpd_{{\rm ad}
\hat{\Lambda}}{\bf R}$ of  Lie algebras  
 ${\goth X}_{H}(M)\cpd_{\beta}{\bf R}$ and ${\bf R}$.

\section{Infinitesimal Takesaki duality}

For an infinitesimal symplectic dynamical system
 $(M,{\bf R},\beta)$, we show the structure
 of second crossed product 
$({\goth X}_{H}(M)\cpd_{\beta}{\bf R})\cpd_{\hat{\beta}}{\bf R}$.
 In order to  use the duality of abelian group ${\bf R}=\hat{{\bf R}}$,
 we restrict the representation space to 
$C^{\infty}(M)\otimes 
{\cal S}({\bf R})\otimes {\cal S}({\bf R})$
 where ${\cal S}({\bf R})\subset C^{\infty}({\bf R})$ is the set of 
Schwartz class functions over ${\bf R}$
 that is, they  satisfy
  $\lim_{|t|\to \infty}|t|^{n-1}|\partial^{m-1}f(t)|=0$ for
 any $n,m\in{\bf N}$.

Recall the Heisenberg Lie algebra $H_{1}$
 is a Lie algebra with a couple of generators
$a,a^{\sdag}$ which satisfying $[a,a^{\sdag}]=I$
 over ${\bf C}$.

By these preparations, we can state our main theorem.

\begin{Thm}(Infinitesimal Takesaki duality)
For an infinitesimal symplectic dynamical system
 $(M,{\bf R},\beta)$ induced by a symplectic dynamical system
 $(M, {\bf R},\alpha)$  
\begin{enumerate}
\item there is a derivative action 
$\delta :H_{1}\to {\rm End}({\goth X}_{H}(M))$
 such that 
\[ \delta_{a}X_{f}        =  \frac{1}{\sqrt{2}}X_{\beta_{1}f},\hspace{.2in}
   \delta_{a^{\sdag}}X_{f}= -\frac{1}{\sqrt{2}}X_{\beta_{1}f}
\hspace{.2in}\mbox{ for } f\in C^{\infty}(M).\]
\item For the derivative action  in (i),
 there is an isomorphism  of Lie algebras
\[ ({\goth X}_{H}(M)\cpd_{\beta}{\bf R})\cpd_{\hat{\beta}}{\bf R}
\cong {\goth X}_{H}(M)\cpd_{\delta}H_{1}\]
 where the second crossed product at left hand side
 is represented on $C^{\infty}(M)\otimes 
{\cal S}({\bf R})\otimes {\cal S}({\bf R})$.
\end{enumerate}
\end{Thm}
{\it Proof}
(i) The existence of $\delta$ is showed 
 in the proof of (ii).

(ii) We prove this statement by 
tracing  the proof of Takesaki duality theorem of von Neumann algebra
 ( vol. II, theorem 13.2.9 in \cite{K-R}).

Let \[ {\cal L}_{1}\equiv 
{\goth X}_{H}(M)\cpd_{\beta}{\bf R},\hspace{.2in} 
{\cal L}_{2}\equiv 
({\goth X}_{H}(M)\cpd_{\beta}{\bf R})\cpd_{\hat{\beta}}{\bf R}\]
\[V_{0}\equiv C^{\infty}(M)\otimes{\cal S}({\bf R}),\hspace{.2in}
 V_{3}\equiv C^{\infty}(M)\otimes 
{\cal S}({\bf R})\otimes {\cal S}({\bf R})\]

By definition, ${\cal L}_{1}$
  is generated by 
\[ \pi(f) \hspace{.5in} \Lambda_{s} 
\hspace{1in} f\in C^{\infty}(M),s\in {\bf R}\]
Let $U_{1} : V_{0}\to V_{0}$ be 
 $(U_{1}\Psi)(p,s^{'})\equiv\Psi(\alpha_{-s^{'}}(p), s^{'})$
 for $\Psi\in V_{0}$. Then $U_{1}$ becomes a linear automorphism
 on $V_{0}$ and 
$U_{1}\pi(f)U_{1}^{-1}=-X_{f}\otimes I$,
 $U_{1}\Lambda_{s}U_{1}^{-1}=\beta_{s}\otimes I+I\otimes l_{s}$.
 where $(l_{s}\chi)(p,s^{'})\equiv \frac{d}{dx}\chi(p,s^{'}-xs)|_{x=0}$
 for $\chi\in {\cal S}({\bf R})$. Because ${\rm Ad}U_{1}$
 preserve Lie bracket, 
${\cal L}_{1}^{'}\equiv ({\rm Ad}U_{1})({\cal L}_{1})\cong{\cal L}_{1}$ as
 a  Lie algebra. We obtain a new Lie algebra ${\cal L}_{1}^{'}$ on $V_{0}$
 with generators
\[ \hspace{.2in}X_{f}\otimes I\hspace{.5in}\beta_{s}\otimes I+I\otimes l_{s}
\hspace{1in} f\in C^{\infty}(M),s\in {\bf R}\]

Next step, 
define 
 $U_{2}: V_{3}\to V_{3}$ by 
$(U_{2}\Phi)(p,s,t)\equiv 
 (W_{t}\otimes I\Phi)(p,s,t)$. Then $U_{2}$
  becomes a linear automorphism on $V_{3}$ and 
$U_{2}\hat{\pi}(R)U_{2}^{-1}=R\otimes I$,
 $U_{2}\hat{\Lambda}_{t}U_{2}^{-1}=
\hat{\Lambda}_{t}+I\otimes w_{t}\otimes I
 =I\otimes I\otimes l_{t}+I\otimes w_{t}\otimes I$ 
  for $R\in {\cal L}_{1}$, $t\in{\bf R}$. So 
 a Lie algebra 
${\cal L}_{3}\equiv 
({\rm Ad}U_{2})({\rm Ad}(U_{1}\otimes I))({\cal L}_{2})$ 
with generators
\[ X_{f}\otimes I\otimes I\hspace{.2in}
\beta_{s}\otimes I\otimes I+I\otimes l_{s}\otimes I
\hspace{.2in}
I\otimes I\otimes l_{t}-I\otimes w_{t}\otimes I
\] 
 for $f\in C^{\infty}(M),t,s,\in {\bf R}$
 is isomorphic to ${\cal L}_{2}$ as a Lie algebra on  $V_{3}$.

Third step, by the fact that the  Fourier transformation 
$T :{\cal S}({\bf R})\to{\cal S}({\bf R})$ 
\[ (Tf)(t)\equiv \frac{1}{\sqrt{2\pi}}\int_{{\bf R}}
e^{\sqrt{-1} ts}f(s)ds\hspace{.3in} f\in {\cal S}({\bf R})\]
is bijection
 and $Tl_{t}T^{-1}=w_{t}$,  a
 Lie algebra  ${\cal L}_{4}\equiv (I\otimes I\otimes T){\cal L}_{3}
(I\otimes I\otimes T^{-1})$ with generators
 \[ X_{f}\otimes I\otimes I\hspace{.2in}
\beta_{s}\otimes I\otimes I+ I\otimes l_{s}\otimes I
 \hspace{.2in}
I\otimes I\otimes w_{t}+I\otimes w_{t}\otimes I\]
 for 
 $f\in C^{\infty}(M),t,s,\in {\bf R}$ 
 is  isomorphic to ${\cal L}_{2}$ as Lie algebra on  $V_{3}$.

Forth step, $U_{3}: V_{3}\to V_{3}$ is defined by
 $(U_{3}\Phi)(p,s,t)\equiv\Phi(p,s-t,t)$ for $\Phi \in V_{3}$,
 $p\in M, t,s\in{\bf R}$. Then $U_{3}$ 
 becomes a linear automorphism on $V_{3}$ and 
 $U_{3}(I\otimes w_{t}\otimes I
+I\otimes I\otimes w_{t})U_{3}^{-1}
=I\otimes w_{t}\otimes I$
 and  other generators are invariant under
 ${\rm Ad}(U_{3})$. 
So,
 ${\cal L}_{4}\equiv U_{3}{\cal L}_{3}U_{3}^{-1}$ with generators
 \[ X_{f}\otimes I\otimes I\hspace{.2in}
\beta_{s}\otimes I\otimes I+I\otimes l_{s}\otimes I
\hspace{.2in}
I\otimes w_{t}\otimes I\]
$
 f\in C^{\infty}(M),t,s,\in {\bf R}$
is isomorphic to ${\cal L}_{2}$. So, naturally
 it is isomorphic to ${\cal L}_{5}$ on $V_{0}$
 with generators
 \[ X_{f}\otimes I\hspace{.2in}
I\otimes l_{s}+\beta_{s}\otimes I \hspace{.2in}
I\otimes w_{t}\hspace{.1in}
 f\in C^{\infty}(M),t,s,\in {\bf R}\]
Then, 
${\cal L}_{6}\equiv U_{1}^{-1}{\cal L}_{5}U_{1}$ on $V_{0}$ 
 has generators 
\[ \pi(f)\hspace{.2in}
I\otimes l_{s}\hspace{.2in}
I\otimes w_{t}\hspace{.1in}
 f\in C^{\infty}(M),t,s,\in {\bf R}\]
 and  it is isomorphic to 
${\cal L}_{2}=
({\goth X}_{H}\cpd_{\beta}{\bf R})\cpd_{\hat{\beta}}{\bf R}$.

By the way, $[w_{t}, l_{s}]=-\sqrt{-1}stI$ on
 ${\cal S}({\bf R})$ for $s,t\in {\bf R}$.
 ( In quantum mechanics, they are  usually denoted by
$w_{t}=\sqrt{-1}t\hat{x}$ and $l_{s}=-s\frac{d}{dx}$ on
 ${\cal S}({\bf R})$.) And 
$[I\otimes w_{t}, \pi(f)]=0,
[I\otimes l_{t}$, $\pi(f)]=[\Lambda_{t}, \pi(f)]=\pi
(\beta_{t}(f))$. 
 So, Lie algebra generated by 
$\{I\otimes l_{s},
I\otimes w_{t}: t,s,\in {\bf R}\}$
becomes   Heisenberg Lie algebra
 $H_{1}$ on $V_{0}$.

Here
 let the Sch\"{o}dinger representation of canonical commutation relation
 on 1-dimensional harmonic oscillator
\[a\equiv \frac{1}{\sqrt{2}}I
\otimes (-\sqrt{-1}w_{1}-l_{1})\hspace{.6in}
a^{\sdag}
\equiv \frac{1}{\sqrt{2}}I
\otimes (-\sqrt{-1}w_{1}+l_{1}).\]
 Then
 $[a,a^{\sdag}]=I$ and
\[[a, \pi(f)]=-\frac{1}{\sqrt{2}}\pi(\beta_{1}(f))\hspace{.5in}
[a^{\sdag}, \pi(f)]=\frac{1}{\sqrt{2}}\pi(\beta_{1}(f))\]
 for $f\in C^{\infty}(M)$. So, let 
$\delta :H_{1}\to {\rm End}(\pi(C^{\infty}(M))$
 be 
$\delta_{b}={\rm ad}(b)$
 for  $b\in H_{1}$,  then $\delta$ becomes derivative action.
 For this derivative action, the semi direct product 
$\pi(C^{\infty}(M))\cpd_{\delta} H _{1}$ is a Lie algebra
 with generators
\[ \pi(f) \hspace{.4in} I\otimes w_{t}\hspace{.4in} I\otimes l_{s}\]
 it is just ${\cal L}_{6}$. 
 By identification, 
$\pi(C^{\infty}(M))\cong C^{\infty}(M)/{\bf C}1_{M}
\cong {\goth X}_{H}(M)$ $\pi(f)\mapsto -X_{f}$, 
 statements follow. \qedh

	\section{Discussion}

We obtain a corollary of theorem 6.1.
\begin{cor}
For a smooth symplectic flow $\alpha: {\bf R}\to {\rm Symp}(M)$
 on a symplectic manifold $M$, there is a derivative action
$\delta:H_{1}\to {\rm End}({\goth  X}_{H}(M))$
 such that 
\[ 
\delta_{a}X_{f}        
=\frac{1}{\sqrt{2}}\frac{d}{dt}X_{\alpha_{t}^{*}f}|_{t=0} \hspace{.2in}
\delta_{a^{\sdag}}X_{f}=
-\frac{1}{\sqrt{2}}\frac{d}{dt}X_{\alpha_{t}^{*}f}|_{t=0}\]
 where $H_{1}$ is a Heisenberg algebra with a couple of generators
$a,a^{\sdag}$ satisfying Heisenberg commutation relation
 and ${\goth X}_{H}(M)$ is  a Lie algebra  of all
 Hamiltonian vector fields $X_{f}$, $f\in C^{\infty}(M)$ on $M$.
\end{cor}

Difference between original Takesaki duality
 and infinite case is follows :

$({\cal R}\cpd_{\alpha}{\bf R})\cpd_{\hat{\alpha}}{\bf R}$
 is separated to ${\cal R}$ and ${\cal L}(L_{2}({\bf R}))$
 by tensor product.  On the other hand,  infinitesimal case
 vector fields and Heisenberg algebra are 
  jointed by derivative action. So,  
  it may be  considered that in the case of original 
Takesaki duality, jointed part is 
 absorbed by strong operator closure of algebraic tensor product  
${\cal R}\otimes {\cal L}(L_{2}({\bf R}))$.

\vv

$\nears${\Large {\bf Appendix }}
\def\labelenumi{\thesection.}
\def\theenumi{Appendix\Alph{section}}
\renewcommand{\thesection}{Appendix\Alph{section}}
\appendix

\setcounter{section}{0}



\section{W$^{*}$-dynamical system and Takesaki duality}\label{section:A}

Let ${\cal R}$ be a von Neumann algebra acting on a Hilbert space ${\cal H}$.
 A W$^{*}$-dynamical system $({\cal R}, G, \alpha)$ is a data such that 
 $G$ is a locally compact group and $\alpha : G\to {\rm Aut}({\cal R})$
 is a strong  operator continuous action of $G$. 

The crossed product ${\cal R}\cpd_{\alpha}G$ of $({\cal R}, G, \alpha)$
 is a von Neumann algebra  acting on a
 Hilbert space  ${\cal H}\otimes L_{2}(G, \mu)\cong
L_{2}(G,{\cal H})$ 
 which has generator $\pi(R), \lambda_{g}$ defined by
\[(\pi(R)\Psi)(h)\equiv((\alpha_{h^{-1}}(R)\otimes I)\Psi)(h)
\hspace{.4in} \lambda_{g}\equiv I\otimes L_{g}\]
 for $g, h\in G$, $R\in {\cal R}$, 
$\Psi\in  {\cal H}\otimes L_{2}(G, \mu)$ 
where $\mu$ is the left Haar measure of $G$, $L_{g}$
 is the left regular representation of $G$.

Assume $G$ is abelian and let $\hat{G}$
 be a dual of $G$.  
 Then we can define the second crossed product
 $({\cal R}\cpd_{\alpha}G)\cpd_{\hat{\alpha}}\hat{G}$ of 
$({\cal R}, G, \alpha)$. 
$({\cal R}\cpd_{\alpha}G)\cpd_{\hat{\alpha}}\hat{G}$
 is defined by a von Neumann algebra 
 acting on a Hilbert space
 ${\cal H}\otimes L_{2}(G,\mu)\otimes
 L_{2}(\hat{G},\nu)$ generated by
 $\hat{\pi}(S)$, $\hat{\lambda}_{\chi}$ 
for $S\in {\cal R}\cpd_{\alpha}G$,
 $\chi\in \hat{G}$  which defined by
\[ (\hat{\pi}(S)\Phi)(h,\eta)\equiv 
((\hat{\alpha}_{\eta^{-1}}(S)\otimes I)\Phi)(h,\eta)
\hspace{.3in}
\hat{\lambda}_{\chi}\equiv I\otimes I\otimes L_{\chi}\]
 for $S\in {\cal R}\cpd_{\alpha}G$, $\chi,\eta\in \hat{G}$,
 $\Phi\in  {\cal H}\otimes L_{2}(G,\mu)\otimes
 L_{2}(\hat{G},\nu)$ 
 where $\nu$ is the left Haar measure of $\hat{G}$,
$L_{\chi}$
 is the left regular representation of $\hat{G}$.

 Then 
\begin{Thm}(Takesaki duality \cite{Takesaki})
 \[({\cal R}\cpd_{\alpha}G)\cpd_{\hat{\alpha}}\hat{G}
\cong {\cal R}\bar{\otimes }{\cal L}(  L_{2}(G,\mu) )\]
where ${\cal L}(L_{2}(G,\mu))$ is 
the algebra of all bounded linear operators on 
$L_{2}(G,\mu)$.
 \end{Thm}
An application of this theorem 
 is  the structure theorem of type III-factor.






\vv

\newpage

\end {document}